\documentclass[12pt]{article}
\usepackage{latexsym,amssymb,amsmath}
\newtheorem{theorem}{Theorem}
  \newtheorem{proposition}{Proposition}

\newcommand{\N}{\mathbb{N}}
\newcommand{\mDelta}{{\mathbf{\Delta}}}
\newcommand{\X}{\mathbb{X}}
\newcommand{\mK}{\mbox{K}}
\newcommand{\K}{\mathbb{K}}

\newcommand{\D}{\mbox{Dim}}

\newcommand{\cH}{{\cal H}}
\newcommand{\R}{{\mathbb{R}}}
\newcommand{\mH}{{\mathbb{H}}}
\newcommand{\C}{{\mathbb C}}
\newcommand{\ua}{\underline{a}}

\newcommand{\ux}{\underline{x}}
\newcommand{\uy}{\underline{y}}

\newcommand{\upa}{\underline{\partial}}
\newcommand{\E}{{\bf E}}
\def\pf{\par\noindent {\em Proof.}~\par\noindent}
\def\qed{~\hfill{$\square$}\pagebreak[1]\par\medskip\par}

\begin{document}


\title{Complex and Quaternionic Cauchy formulas in Koch snowflakes}
\author
{Marisel Avila Alfaro$^{(1)}$; Ricardo Abreu Blaya$^{(2)}$}
\vskip 1truecm
\date{\small $^{(1)}$ Instituto de Matem\'aticas, Universidad Nacional Aut\'onoma de M\'exico, M\'exico\\$^{(2)}$ Facultad de Matem\'aticas, Universidad Aut\'onoma de Guerrero, M\'exico.\\Emails: mavilaa94@comunidad.unam.mx, rabreublaya@yahoo.es}
\maketitle
\begin{abstract}
In this paper we derive a Cauchy integral formula for holomorphic and hyperholomorphic functions in domains bounded by a Koch snowflake in two and three dimensional setting. 
\end{abstract}
\vspace{0.3cm}

\small{
}

\section{Introduction}
The Cauchy integral formula plays a decisive  role in Complex Analysis. An important consequence shows us a surprising property of holomorphic functions, namely that these are real continuously differentiable not only once, but infinitely many times. Thus all derivatives of a holomorphic function are again holomorphic.

Cauchy  formula says that a holomorphic function $f$ is completely defined by its boundary values. More precisely, let $\Omega$ be a Jordan domain in $\R^2$ with sufficiently smooth boundary $\Gamma$ and let the function $f$ be holomorphic in $\Omega$ and continuous in $\overline{\Omega}$, then
\begin{equation}\label{CCF}
\frac{1}{2\pi i}\int\limits_{\Gamma}\frac{f(\xi)}{\xi-z}d\xi=\biggl\{
\begin{array}{rl}
f(z),\,\,& z\in\Omega_+\\
0,\,\,& z\in\Omega_-,
\end{array} 
\end{equation} 
where here and in the sequel use has been made of the notation $\Omega_+=\Omega$, $\Omega_-=\overline{\Omega}^c$.

Wider scope could be obtained when we replace complex-valued functions by functions defined in $\R^3$ and taking their
values in the skew-field of quaternions. In the so-called Quaternionic Analysis, an analogous multidimensional Cauchy formula can be derived as an easy consequence of the Stokes Theorem. In this context this formula, whose detailed notational meaning will be explained below, takes the form 
\begin{equation}\label{QCF}
\int\limits_{\Gamma}E_0(\uy-\ux)\nu(\uy)v(\uy)d\uy=\biggl\{
\begin{array}{rl}
v(z),\,\,& \ux\in\Omega_+\\
0,\,\,& \ux\in\Omega_-,
\end{array} 
\end{equation} 
where $v$ is assumed to be a quaternion-valued hyperholomorphic function in the Jordan domain $\Omega$, this time in $\R^3$.

It should be pointed out that both complex and quaternionic version of the Cauchy formula are classically considered in domains with sufficiently smooth boundary \cite{GS,HGS}. It is then natural to ask whether this important formula remains valid (in some sense) in the case of $\Gamma$ being a fractal. This is the question we shall be concerned in this work. We restrict ourselves to the concrete but important case of Koch snowflakes in two and three-dimensional setting. However, our results are expected to be extended to more general cases of deterministic fractals defined as the fixed point of a so-called Iterated Function System \cite{Br}.
\section{Preliminaries}
We do not discuss here the basic facts of complex numbers and complex-valued functions, because we consider them a very classical matter. Quaternions and quaternion-valued functions take up mostly the whole of this brief preliminary section.

Let $\mH$ be the skew field of real quaternions generated by the basic elements $e_0=1,e_1, e_2, e_3$ that fulfill the condition
\[e_ie_j+e_je_i=-2\delta_{ij},\;i,j=1,2,3\]
\[e_1e_2=e_3;\;e_2e_3=e_1;\;e_3e_1=e_2.\]
For each $a=a_0+\ua$, where $\ua:=\sum_{j=1}^{3}a_je_j$, $a_j\in\R$, the norm of $a$ is defined to be $|a|^2=\sum_{j=0}^3a_j^2$. What is more, for $a, b$ from $\mH$ there hold: $|ab|=|a||b|$. The conjugate element $\bar{a}$ is given by $\bar{a}=a_0-\ua$ and we have the properties $a\overline{a}=\overline{a}a=|a|^2$; $a_0:=Sc\;a$ is called scalar part and $\ua= Vec\,a$ is called vector part of the quaternion $a$. If $Sc\,a=0$ the quaternion $a$ is called pure imaginary one, and it is identified with a vector $\vec{a}$ from $\R^3$.

In vector terms, the multiplication of two arbitrary quaternions $a,b$ can be rewritten as follows
\[ab=a_0b_0-\vec{a}\cdot\vec{b}+a_0\vec{b}+b_0\vec{a}+\vec{a}\times\vec{b},\]
where $\vec{a}\cdot\vec{b}$ denotes the scalar product in $\R^3$ and $\vec{a}\times\vec{b}$ denotes Gibbs cross product in $\R^3$. Quaternion multiplication is associative, distributive, but not commutative.

We will consider functions defined on subsets of $\R^3$ and taking values in $\mH$. Those functions might be written as 
\[u(\ux)=\sum_{j=0}^3u_j(\ux)e_j,\]
where the $u_j$'s are $\R$-valued functions. The notions of continuity, differentiability and integrability of a $\mH$-valued function $u$ have the usual component-wise meaning. In particular, the spaces of $\alpha$-H\"older continuous, $k$-time continuous differentiable and $p$-integrable functions are denoted by $C^{0,\alpha}(\E)$, $C^k(\E)$ and $L_p(\E)$ respectively, where $\E$ is a given subset of $\R^3$. 

The so-called Dirac operator $\upa$ is defined to be
\[
\upa=e_1\frac{\partial}{\partial x_1}+e_2\frac{\partial}{\partial x_2}+e_3\frac{\partial}{\partial x_3}.
\]
An $\mH$-valued function $u$, defined and differentiable in an open region $\Omega\subset\R^3$, is called left-hyperholomorphic (right-hyperholomorphic) if $\upa u=0$ ($u\upa=0$) in $\Omega$. Whenever there is no danger of confusion, we will refer to left hyperholomorphic functions simply as hyperholomorphic. The fundamental solution of the operator $\upa$ is given by
\begin{displaymath}
E_0(\ux)=-\frac{1}{4\pi}\frac{\ux}{|\ux|^3} \,\,\,\,\, (\ux\neq 0),
\end{displaymath}
which is commonly referred to as Cauchy kernel.

A quaternionic reformulation of the Gauss theorem yields a sort of Cauchy integral theorem for hyperholomorphic functions. Let $\Omega\subset\R^m$ with sufficiently smooth boundary $\Gamma=\partial\Omega$ and let $u,v\in C^1(\overline{\Omega})$. If moreover $u$ ($v$) is right hyperholomorphic (left-hyperholomorphic) in $\Omega$, then
\begin{equation}\label{QIT}
\int\limits_\Gamma u(\uy)\nu(\uy)v(\uy)d\uy=0,
\end{equation}
where $\nu(\uy)$ stands for the outer normal vector to $\Omega$ at $\uy\in\Gamma$.

Similarly, by standard  arguments we derive the Borel-Pompeiu representation formula
\begin{equation}\label{BP}
\int\limits_{\Gamma}E_0(\uy-\ux)\nu(\uy)v(\uy)d\uy-\int\limits_{\Omega}E_0(\uy-\ux)\upa v(\uy)d\uy=\biggl\{
\begin{array}{rl}
v(z),\,\,& \ux\in\Omega_+\\
0,\,\,& \ux\in\Omega_-
\end{array} 
\end{equation}
for every $\mH$-valued function $v\in C^1(\overline{\Omega})$.

In particular, for hyperholomorphic functions in $\Omega$, \eqref{BP} yields \eqref{QCF}, as a reminiscence of the classical Cauchy integral formula for holomorphic functions.

Thus, after this rapid review of some of the fundamental facts concerning quaternionic analysis, we turn to a very deep theorem from real analysis due to H. Whitney \cite{Wh}, which represents a crucial technical ingredient in our achievement.
\begin{theorem}[Whitney's Extension Theorem]\label{Wh}
Let $\E\subset\R^m$ be compact and a real-valued function $f \in C^{0,\alpha}(\E)$, $0<\alpha\le 1$. Then there exists a compactly supported function $\tilde{f}\in C^{0,\alpha}(\R^m)$ and satisfying
\begin{itemize}
\item[(i)] $\tilde{f}|_\E=f|_\E$;
\item[(ii)]$\bigg|\dfrac{\partial\tilde{f}(\ux)}{\partial x_j}\bigg|\le c{\mbox{dist}}(\ux,\E)^{\alpha-1}$ for $\ux\in\R^m\setminus\E$.    
\end{itemize} 
\end{theorem}
Here and until the end of this work, $c$ will denote a positive constant, not necessarily the same at different occurrences.

For a thorough treatment of the previous theorem we refer the reader to \cite[Chapter VI]{St}.
\section{Basis of fractal geometry: Koch-type fractals}
Following \cite{Br}, fractals live in the space $\cH(\X)$ of compact subsets of a complete metric space $({\X}, {\bf d})$. For the purpose of this paper it is sufficient to assume that $\X=\R^m$, $m=2,3$ and ${\bf d}$ is the usual Euclidean distance $|\ux-\uy|$. As is well-known the space $\cH(\X)$ becomes a complete metric space with the Hausdorff distance
\begin{equation}\label{Hd}
{\bf h}(A,B)= \max\{{\bf d}(A,B), {\bf d}(B,A)\},
\end{equation}
 where 
\[
{\bf d}(A,B):=\max\{{\bf d}(a,B):\,a\in A\}
\]
for $A,B\in\cH(\X)$. It is to be noted that in general ${\bf d}(A,B)\not={\bf d}(B,A)$.

A classical way of comparing fractals is by means of various numbers, which are generally referred to as fractal dimensions. Although there are several of such numbers to measure the roughness of a fractal, we restrict ourselves to the so-called box dimension \cite{Fa}. Let $A\in\cH(\X)$ and denote by $N(A,\epsilon)$ the smallest number of closed balls of radius $\epsilon$ needed to cover $A$. The quantity   
\[
\D(A):=\limsup_{\epsilon\to 0}\frac{\log N(A,\epsilon)}{-\log\epsilon}
\]
is called the box dimension of $A$. This terminology is justified by the fact that the above limit is unchanged if $N(A,\epsilon)$ is thinking as the number of closed square boxes of side length $2^{-k}$, with $2^{-k}\le \epsilon<2^{-k+1}$, which intersect $A$.
\subsection{Koch Snowflake}
The Koch snowflake, being nowhere differentiable, is an example of continuous closed curve in $\R^2$ whose box dimension is given by $\dfrac{\ln 4}{\ln 3}$. This important structure is not only relevant from a mathematical point of view, but also has important applications in engineering and is widely used in modern telecommunication systems \cite{Be,Kar,Tu}. 

The Koch snowflake is constructed by starting with an equilateral triangle $\mK_0$ of side length $l$. Divide each side of $\mK_0$  into three segments of equal length and then remove the middle third of each line and replace it by two line segments of length $\dfrac{l}{3}$ pointing outward to make a corner. The resulting set is a closed curve $\mK_1$ consisting of $12$ line segments of length $\dfrac{l}{3}$. Then with $\mK_n$ denoting the $n$-th iteration, we see that the length of each side of $\mK_n$ is $\dfrac{l}{3^n}$. Moreover we also see easily that the perimeter of $\mK_n$ is equal to $\dfrac{4^n}{3^{n-1}}l$. Repeating the above procedure infinitely many times, we obtain the so-called Koch snowflake $\mK$. It should be noticed that ${\bf h}(\mK_n,\mK)\to 0$ as $n\to\infty$. That this indeed is satisfied follows from the fact that ${\bf d}(\mK_n,\mK)$ and ${\bf d}(\mK,\mK_n)$ are both dominated by $\dfrac{l}{3^{n}}$. 

A natural three-dimensional extension of the the Koch snowflake is similarly constructed as follows. Start with a  regular tetrahedron $\K_0\subset\R^3$ of edge length $l$ and then build a regular tetrahedron of edge length $\dfrac{l}{2}$ based on the midpoints of the forth triangular faces in $\K_0$. After removing the bases of each one of these tetrahedra, the first iteration $\K_1$ is obtained. The Koch snowflake  $\K$ in $\R^3$ is thus the limit approached as the above iterations are followed indefinitely.  Therefore, denoting by $\K_n$ the $n$-th iteration of $\K$, one obtains ${\bf h}(\K_n,\K)\to 0$ as $n\to\infty$ as was to be expected.  Moreover, the fractal dimension of the Koch snowflake  $\K$ is $\dfrac{\log 6}{\log 2}$.  
\section{Integration and Cauchy formula}
The techniques and arguments used in this section are, in spirit, almost independent of dimension. However, the non-commutativity of the product in $\mH$ represents an additional difficulty and in part for that reason we consider the two and three-dimensional cases separately.    
\subsection{Integration and complex Cauchy formula in $\mK$}
We begin by defining a natural integration of suitable continuous functions defined on the Koch snowflake. Suppose $f(z)$ be a complex-valued function in $C^{0,\alpha}(\mK)$ and let $\tilde{f}$ stands for a Whitney extension of $f$. Since by Theorem \ref{Wh}, $\tilde{f}$ is defined and H\"older continuous in the whole $\R^2$, we are tempted to define by a natural limiting argument the integral of $f$ on $\mK$ as
\begin{equation}\label{IK2}
\int\limits_{\mK} f(z)dz^{*}=\lim\limits_{n\to\infty} \int\limits_{\mK_n}\tilde{f}(z)dz,
\end{equation}
if the limit exists.

To make this hypothetical definition unambiguous, we need to ensure that the above limit will not depend on the Whitney extension $\tilde{f}$, which is not uniquely defined by $f$. The following proposition provides a sufficient condition  to legitimate \eqref{IK2}. 
\begin{proposition}\label{legitimo2}
Let $f\in C^{0,\alpha}(\mK)$, with $\alpha$ satisfying   
\begin{equation}\label{>}
\alpha+1>\dfrac{\log 4}{\log 3}.
\end{equation} 
Then the limit \eqref{IK2} exists and its value does not depend on the choice of $\tilde{f}$.
\end{proposition}
\pf In order to show that the limit \eqref{IK2} exists it suffices to see that $\int\limits_{\mK_n}\tilde{f}(z)dz$ is a Cauchy sequence in $\C$. In fact, for fixed $p\in\N$ we have:
\begin{displaymath}
\int\limits_{K_{n+p}}\tilde{f}(z)dz-\int\limits_{K_n}\tilde{f}(z)dz=\sum\limits_{k=0}^{p-1}[\int\limits_{K_{n+p-k}}\!\!\tilde{f}(z)dz-\!\!\!\int\limits_{K_{n+p-k-1}}\!\tilde{f}(z)dz].
\end{displaymath}
We are reduced to prove that each term in the above sum converges to $0$ as $n\to\infty$. To this end and for simplicity, write $m=n+p-k$.
   
It follows from the construction of $\mK$ that
\[
\int\limits_{K_{m}}\!\!\tilde{f}(z)dz-\!\!\!\int\limits_{K_{m-1}}\!\tilde{f}(z)dz=\sum\limits_{r=1}^{3\cdot 4^{m-1}}\int\limits_{\Delta_{r}}\tilde{f}(z)dz,
\]
where the sets $\Delta_{r}$ are equilateral triangles of side length $\frac{l}{3^m}$.

Let $z_{r}\in\Delta_{r}$ be fixed, then Cauchy integral theorem leads to
\[\int\limits_{\Delta_{r}}\tilde{f}(z_{r})dz=0\] 
for every $r=1,2,\dots,{3\cdot 4^{m-1}}$.

Then
\begin{displaymath}
\vert\sum\limits_{r=1}^{3\cdot 4^{m-1}}\int\limits_{\Delta_{r}}\tilde{f}(z)dz\vert
=\vert\sum\limits_{r=1}^{3\cdot 4^{m-1}}\int\limits_{\Delta_{r}}[\tilde{f}(z)-\tilde{f}(z_{r})]dz\vert 
\leq \sum\limits_{r=1}^{3\cdot 4^{m-1}}\int\limits_{\Delta_{r}}\vert \tilde{f}(z)-\tilde{f}(z_{r})\vert \vert dz\vert
\end{displaymath}
and hence
\begin{displaymath}
\sum\limits_{r=1}^{3\cdot 4^{m-1}}\int\limits_{\Delta_{r}}\vert \tilde{f}(z)-\tilde{f}(z_{r})\vert \vert dz\vert \leq c\sum\limits_{r=1}^{3\cdot 4^{m-1}}\int\limits_{\Delta_{r}}\vert z-z_{r}\vert^{\alpha}\vert dz\vert,
\end{displaymath}
where use has been made of the H\"older continuity of $\tilde{f}$. 
 
Since for $z\in\Delta_{r}$ we have $\vert z-z_{r}\vert \leq \dfrac{l}{3^{m}}$, it follows that  
\begin{displaymath}
\sum\limits_{r=1}^{3\cdot 4^{m-1}}\int\limits_{\Delta_{r}}\vert z-z_{r}\vert ^{\alpha}\vert dz\vert \leq \sum\limits_{r=1}^{3\cdot 4^{m-1}}\dfrac{l^{\alpha}}{3^{m\alpha}} \int\limits_{\Delta_{r}}\vert dz\vert=\dfrac{9l^{\alpha +1}}{4}\bigg(\dfrac{4}{3^{(\alpha +1)}}\bigg)^m,
\end{displaymath}
which obviously converges to $0$ if \eqref{>} is satisfied.

We shall show next that the integral given by \eqref{IK2} does not depend on the Whitney extension $\tilde{f}$. Suppose $\tilde{f}$ and $\tilde{\tilde{f}}$ are two different Whitney extension of $f$ and prove that
\begin{displaymath}
\lim\limits_{n \rightarrow \infty} \int\limits_{K_n}\tilde{\tilde{f}}(z)dz = \lim\limits_{n \rightarrow \infty}   \int\limits_{K_n}\tilde{f}(z)dz,
\end{displaymath}
or equivalently that
\begin{displaymath}
\lim\limits_{n \rightarrow \infty} \int\limits_{K_n}[\tilde{f}(z)- \tilde{\tilde{f}}(z)]dz=0.
\end{displaymath}

Of course the function $\hat{f}(z):=\tilde{f}(z)- \tilde{\tilde{f}}(z)$ is H\"older continuous in $\R^2$ with vanishing restriction to $\mK$, i.e., $\hat{f}|_{\mK}=0 $.

Let $K_{n,m}$, $m=1,\dots,3\cdot 4^{n}$, stand for the sides (line segments) of $K_n$ and denote by $z_{n,m}$ one of its ends. Obviously, $z_{n,m}$ is in $\mK$ and for $z\in K_{n,m}$ we have $|z-z_{n,m}|\le\dfrac{l}{3^m}$.

Therefore the integrals
\begin{displaymath}
\int\limits_{K_n}\hat f(z)dz=\sum\limits_{m=1}^{3\cdot 4^{n}}\int\limits_{K_{n,m}} \hat{f}(z)dz= \sum\limits^{3\cdot 4^{n}}\limits_{m=1} \int\limits_{K_{n,m}}[\hat f(z)-\hat f(z_{n,m})]dz
\end{displaymath}
are dominated in absolute value by 
\begin{displaymath}
\sum\limits^{3\cdot 4^{n}}\limits_{m=1}\int\limits_{K_{n,m}}\vert \hat f(z)-\hat f(z_{n,m})\vert \vert dz\vert\leq c\sum\limits^{3\cdot 4^{n}}\limits_{m=1}\int\limits_{K_{n,m}}\vert z - z_{n,m}\vert^{\alpha}\vert dz\vert. 
\end{displaymath}
Arguing in the same way as in the proof of the existence, one finds that under condition \eqref{>} the last term tends to zero as $n\to\infty$, and we are done.\qed

We let $\mK_+$ ($\mK_-$) to be the bounded (unbounded) open domain determined by $\mK$. Let us now state and prove a sort of Cauchy integral theorem in our context.
\begin{theorem}\label{CCT}
Let $f\in\ C^{0,\alpha}(\mK_+\cup\mK)$ with $\alpha$ satisfying \eqref{>}. If moreover $f$ is holomorphic in $\mK_+$, then
\begin{equation}\label{CIT}
\int\limits_{\mK}f(z)dz^{*}=0
\end{equation}
\end{theorem}
\pf By definition we have
\[
\int\limits_{\mK}f(z)dz^{*}=\lim\limits_{n\to\infty} \int\limits_{\mK_n}[\tilde{f}(z)-f(z)]dz,
\] 
because $\int\limits_{\mK_n}f(z)dz=0$ by the usual Cauchy theorem.

Since $\tilde{f}(z)-f(z)$ belongs to $C^{0,\alpha}(\mK_+\cup\mK)$ and vanishes at all points of $\mK$, a quite similar arguments to those used in the proof of Proposition \ref{legitimo2} yields the desired conclusion.\qed

Next, we have in mind a Cauchy integral formula formulated as follows.
\begin{theorem}\label{*TCF}
Let $f\in\ C^{0,\alpha}(\mK_+\cup\mK)$ with $\alpha$ satisfying \eqref{>}. If moreover $f$ is holomorphic in $\mK_+$, then
\begin{equation}\label{*CCF}
\frac{1}{2\pi i}\int\limits_{\mK}\frac{f(\xi)}{\xi-z}d\xi^{*}=\biggl\{
\begin{array}{rl}
f(z),\,\,& z\in\mK_+\\
0,\,\,& z\in\mK_-,
\end{array} 
\end{equation} 
\end{theorem} 
\pf Suppose $z\in\mK_+$ and let $d(z,\mK)$ be the distance of $z$ to $\mK$. Then for $\xi_1,\xi_2\in\mK$
\[
|\frac{1}{\xi_1-z}-\frac{1}{\xi_2-z}|\le\frac{|\xi_1-\xi_2|}{|\xi_1-z||\xi_2-z|}\le d(z,\mK)^{-2}|\xi_1-\xi_2|.
\]
Hence, the function $\phi(\xi):=\frac{f(\xi)}{\xi-z}$ is obviously in $C^{0,\alpha}(\mK)$. Then by \eqref{IK2}, whose application is again legitimate in view of the above property of $\phi$, we have
\[
\frac{1}{2\pi i}\int\limits_{\mK}\frac{f(\xi)}{\xi-z}d\xi^{*}=\lim\limits_{n\to\infty}\frac{1}{2\pi i}\int\limits_{\mK_n}\tilde{\phi}(\xi)d\xi
\]
Choose $N$ large enough so that for $n\ge N$ the point $z$ lies in the bounded open domain determined by the Jordan curve $\mK_N$. For such $n$ we have
\[
\frac{1}{2\pi i}\int\limits_{\mK_n}\tilde{\phi}(\xi)d\xi=\frac{1}{2\pi i}\int\limits_{\mK_n}[\tilde{\phi}(\xi)-\frac{{f}(\xi)}{\xi-z}]d\xi+\frac{1}{2\pi i}\int\limits_{\mK_n}\frac{{f}(\xi)}{\xi-z}d\xi
\]
or equivalently
\begin{equation}\label{equiv}
\frac{1}{2\pi i}\int\limits_{\mK_n}\tilde{\phi}(\xi)d\xi=\frac{1}{2\pi i}\int\limits_{\mK_n}[\tilde{\phi}(\xi)-\frac{{f}(\xi)}{\xi-z}]d\xi+f(z),
\end{equation}
where we have used the usual Cauchy integral formula applied to $f$ in $\mK_n$.

Moreover we notice that the function $\frac{1}{\xi-z}$ is uniformly Lipschitz continuous in $\mK_n$, $n\ge N$, with uniform Lipschitz constant $d(z,\mK_N)^{-2}$. This fact is a direct consequence of the obvious relation
\[
d(z,\mK_N)\le d(z,\mK_n),\,\mbox{for}\,\,n\ge N,
\]
which is itself implied by the construction of $\mK$.

Consequently  $\tilde{\phi}(\xi)-\phi(\xi)$ belongs (uniformly) to $C^{0,\alpha}(\mK_n)$,  for $n\ge N$, and its restriction to $\mK$ vanishes. With this in mind we can apply the reasoning already used in Proposition \ref{legitimo2} to get
\[
\lim\limits_{n\to\infty}\frac{1}{2\pi i}\int\limits_{\mK_n}[\tilde{\phi}(\xi)-\frac{{f}(\xi)}{\xi-z}]d\xi=0,
\]
which leads to \eqref{*CCF}.

For $z\in\mK_-$ matters become obvious by Theorem \ref{CCT}.\qed  
\subsection{Integration and quaternionic Cauchy formula in $\K$}
Here we introduce an integration of $\mH$-valued  H\"older continuous functions over the snowflake $\K$ in $\R^3$, which is constructed along the same lines we have used in the previous subsection.

As before, the construction of $\K$ and the Whitney extension theorem (applied componentwise) suggest the use of a limiting argument yielding a tentative integral of $u$ on $\K$ given by
\[
\int\limits_{\K} u(\uy)d\uy^*=\lim\limits_{n\to\infty}\int\limits_{\K_n}\tilde{u}(\uy)\nu_n(\uy)d\uy,
\]
if the limit exists. Here and below $\nu_n(\uy)$ stands for the outer normal vector at $\uy\in\K_n$. Since $\K_n$ is suffciently smooth (Lipschitz, for instance), such a vector $\nu_n(\uy)$ exists almost everywhere in $\K_n$.

Alternatively, due tho the non-commutativity of the product in $\mH$, one could consider the integral on $\K$ as 
\[
\int\limits_{\K}d\uy^* u(\uy)=\lim\limits_{n\to\infty}\int\limits_{\K_n}\nu_n(\uy)\tilde{u}(\uy)d\uy.
\]
However, a closer look to quaternionic integration theory reveals that it usually involves two functions at once, one on the left and another on the right of the outer normal vector $\nu$. This gives rise to a more appropriate definition given by
\begin{equation}\label{IK3}
\int\limits_{\K} u(\uy)d\uy^* v(\uy)=\lim\limits_{n\to\infty}\int\limits_{\K_n}\tilde{u}(\uy)\nu_n(\uy)\tilde{v}(\uy)d\uy,
\end{equation}
where $u,v\in C^{0,\alpha}(\K)$.

The proof of the following statement is similar to that of Proposition \ref{legitimo2}, for that reason it will be, for the most part, only sketched.
\begin{proposition}\label{legitimo3}
Let $u,v\in C^{0,\alpha}(\K)$, with $\alpha$ satisfying   
\begin{equation}\label{>H}
\alpha+2>\dfrac{\log 6}{\log 2}.
\end{equation} 
Then the limit \eqref{IK3} exists and its value does not depend on the choice of $\tilde{u}$, $\tilde{v}$.
\end{proposition} 
\pf We begin by showing that
\[
\int\limits_{\K_n}\tilde{u}(\uy)\nu_n(\uy)\tilde{v}(\uy)d\uy
\]
is a Cauchy sequence in $\mH$. For this purpose we need to prove that the sum
\[
\sum\limits_{r=1}^{4\cdot 6^{m-1}}\int\limits_{\mDelta_{r}}\tilde{u}(\underline{y})\nu(\uy)\tilde{v}(\uy)d\uy
\]
tends to zero as $m\to\infty$, where $\mDelta_{r}$ denotes a regular tetrahedron of edge length $\dfrac{l}{2^m}$. Indeed choosing a fixed point $\uy_r\in\mDelta_r$ we have
\begin{eqnarray}\label{importante}
\sum\limits_{r=1}^{4\cdot 6^{m-1}}\int\limits_{\mDelta_{r}}\tilde{u}(\uy)\nu(\uy)\tilde{v}(\uy)d\uy=\sum\limits_{r=1}^{4\cdot 6^{m-1}}\int\limits_{\mDelta_{r}}[\tilde{u}(\underline{y})-\tilde{u}(\uy_r)]\nu(\uy)\tilde{v}(\uy)d\uy\nonumber\\+\sum\limits_{r=1}^{4\cdot 6^{m-1}}\int\limits_{\mDelta_{r}}\tilde{u}(\uy_r)\nu(\uy)[\tilde{v}(\underline{y})-\tilde{v}(\uy_r)]d\uy,
\end{eqnarray}
where we used the fact that
\[
\int\limits_{\mDelta_{r}}\tilde{u}(\underline{y}_r)\nu(\uy)\tilde{v}(\underline{y}_r)d\uy=0
\]
as a rather trivial implication of \eqref{QIT}.

Consequently it follows from \eqref{importante} and from the fact that $\tilde{u},\tilde{v}$ are bounded, that
\begin{eqnarray*}
|\sum\limits_{r=1}^{4\cdot 6^{m-1}}\int\limits_{\mDelta_{r}}\tilde{u}(\uy)\nu(\uy)\tilde{v}(\uy)d\uy|\le c_1\sum\limits_{r=1}^{4\cdot 6^{m-1}}\int\limits_{\mDelta_{r}}|\tilde{u}(\uy)-\tilde{u}(\uy_r)|d\uy\\+c_2\sum\limits_{r=1}^{4\cdot 6^{m-1}}\int\limits_{\mDelta_{r}}|\tilde{v}(\underline{y})-\tilde{v}(\uy_r)|d\uy\le c\,\dfrac{l^{\alpha + 2}\cdot\sqrt{3}\cdot 6^{m}}{2^{m(\alpha+2)}},
\end{eqnarray*}
which under condition  \eqref{>H} tends to zero as $m\to\infty$.

Next we turn to the independence of the Whitney extension. Suppose $\tilde{u},\tilde{\tilde{u}}$ and  $\tilde{v},\tilde{\tilde{v}}$ are two different Whitney extensions of $u$ and $v$, respectively. We are reduced to prove that
\[
\lim\limits_{n\to\infty}\int\limits_{\K_n}[\tilde{u}(\uy)\nu_n(\uy)\tilde{v}(\uy)-\tilde{\tilde{u}}(\uy)\nu_n(\uy)\tilde{\tilde{v}}(\uy)]d\uy=0.
\]
In fact, we have 
\begin{eqnarray}\label{truco}
\int\limits_{\K_n}[\tilde{u}(\uy)\nu_n(\uy)\tilde{v}(\uy)-\tilde{\tilde{u}}(\uy)\nu_n(\uy)\tilde{\tilde{v}}(\uy)]d\uy=\int\limits_{\K_n}[\tilde{u}(\uy)-\tilde{\tilde{u}}(\uy)]\nu_n(\uy)\tilde{v}(\uy)d\uy\nonumber\\+\int\limits_{\K_n}\tilde{\tilde{u}}(\uy)\nu_n(\uy)[\tilde{v}(\uy)-\tilde{\tilde{v}}(\uy)]d\uy.
\end{eqnarray}
Of course, $\tilde{u}(\uy)-\tilde{\tilde{u}}(\uy)$ belongs to $C^{0,\alpha}(\R^3)$ and vanishes at all points of $\K$. This together with the boundedness of $\tilde{v}$ allows us to proceed in analogy to the complex case (Proposition \ref{legitimo2}), where the $3\cdot4^n$ line segments are this time replaced by the $4\cdot 6^n$ triangular faces of $\K_n$. With this in mind we easily deduce that the first summand in the right-hand side of \eqref{truco} tends to zero as $n\to\infty$. A quite similar argument gives the same result for the second one.\qed 

In the sequel it will be assumed that $\alpha$ satisfies the condition \eqref{>H}. Moreover, in the customary notation used earlier in the complex case, we will denote by $\K_+$ ($\K_-$)  the bounded (unbounded) open domain in $\R^3$ determined by $\K$.
\begin{theorem}
Let $u,v\in C^{0,\alpha}(\K_+\cup\K)$ be right- (respectively left-) hyperholomorphic $\mH$-valued functions in $\K_+$. Then
\[
\int\limits_{\K} u(\uy)d\uy^* v(\uy)=0
\]
\end{theorem}
\pf We are reduced to prove that
\[
\lim\limits_{n\to\infty}\int\limits_{\K_n} \tilde{u}(\uy)\nu_n(\uy)\tilde{v}(\uy)d\uy=0.
\]
We have
\begin{eqnarray*}
\int\limits_{\K_n} \tilde{u}(\uy)\nu_n(\uy)\tilde{v}(\uy)d\uy=\int\limits_{\K_n}[\tilde{u}(\uy)-u(\uy)]\nu_n(\uy)\tilde{v}(\uy)d\uy\\+\int\limits_{\K_n} u(\uy)\nu_n(\uy)[\tilde{v}(\uy)-v(\uy)]d\uy+\int\limits_{\K_n} u(\uy)\nu_n(\uy)v(\uy)d\uy
\end{eqnarray*}
or equivalently
\begin{eqnarray}\label{truco2}
\int\limits_{\K_n} \tilde{u}(\uy)\nu_n(\uy)\tilde{v}(\uy)d\uy=\int\limits_{\K_n}[\tilde{u}(\uy)-u(\uy)]\nu_n(\uy)\tilde{v}(\uy)d\uy\nonumber\\+\int\limits_{\K_n} u(\uy)\nu_n(\uy)[\tilde{v}(\uy)-v(\uy)]d\uy,
\end{eqnarray}
since 
\[
\int\limits_{\K_n} u(\uy)\nu_n(\uy)v(\uy)d\uy=0
\]
by \eqref{QIT}.
Thus the only thing for us to prove here is that both integrals in the right side of \eqref{truco2} go to zero as $n\to\infty$. The proof of that is merely a repetition of that of  the Proposition \ref{legitimo3}.\qed

Let us now consider a quaternionic Cauchy formula in our fractal context. Technically, the arguments we use in deriving such a formula are quite  similar to those employed in Theorem \ref{*TCF}, but we decided to include its proof here for the sake of completeness.
\begin{theorem}
Let $v\in C^{0,\alpha}(\K_+\cup\K)$ be left-hyperholomorphic in $\K_+$. Then
\begin{equation}\label{QCF*}
\int\limits_{\K}E_0(\uy-\ux)d\uy^* v(\uy)=\biggl\{
\begin{array}{rl}
v(z),\,\,& \ux\in\K_+\\
0,\,\,& \ux\in\K_-,
\end{array} 
\end{equation}
\end{theorem}
\pf Suppose $\ux\in\K_+$ to be fixed and let $d(\ux,\K)$ be the distance between $\ux$ and $\K$. Because of
\[
|E_0(\uy_1-\ux)-E_0(\uy_2-\ux)|\le c|\uy_1-\uy_2|\sum_{i=1}^{2}\frac{1}{|\ux-\uy_1|^i|\ux-\uy_2|^{3-i}}
\]
it follows that
\[
|E_0(\uy_1-\ux)-E_0(\uy_2-\ux)|\le\frac{2c}{d(\ux,\K)^3}|\uy_1-\uy_2|
\]
for $\uy_1,\uy_2\in\K$.

Then the function $\psi(\uy)=E_0(\uy-\ux)$ belongs to $C^{0,\alpha}(\K)$ and by definition
\[
\int\limits_{\K}E_0(\uy-\ux)d\uy^* v(\uy)=\lim\limits_{n\to\infty}\int\limits_{\K_n}\tilde{\psi}(\uy)\nu_n(\uy)\tilde{v}(\uy)d\uy.
\]
As before choose $N$ large enough so that for $n\ge N$ the point $\ux$ lies in the bounded open domain determined by $\K_N$. 

For such $n$ we have
\begin{eqnarray}\label{truco3}
\int\limits_{\K_n}\tilde{u}(\uy)\nu_n(\uy)\tilde{v}(\uy)d\uy=\int\limits_{\K_n}[\tilde{\psi}(\uy)-E_0(\uy-\ux)]\nu_n(\uy)\tilde{v}(\uy)d\uy+\int\limits_{\K_n}E_0(\uy-\ux)\nu_n(\uy)\tilde{v}(\uy)d\uy\nonumber\\=\int\limits_{\K_n}[\tilde{\psi}(\uy)-E_0(\uy-\ux)]\nu_n(\uy)\tilde{v}(\uy)d\uy+\int\limits_{\K_n}E_0(\uy-\ux)\nu_n(\uy)[\tilde{v}(\uy)-v(\uy)]d\uy+v(\ux),
\end{eqnarray}
the last equality being a direct consequence of the quaternionic Cauchy formula \eqref{QCF}.

We are therefore faced with the task of showing that the last two integrals in \eqref{truco3} tend to zero as $n\to\infty$. For this it will be sufficient to repeat the reasoning leading to \eqref{*CCF}, but using the fact that both $\tilde{\psi}-\psi$ and $\tilde{v}(\uy)-v(\uy)$ belong (uniformly) to $C^{0,\alpha}(\K_n)$, $n\ge N$, and their restrictions to $\K$ vanish.

For $\ux\in\K_-$ the proof is again trivial and will be omitted.\qed
\section{Concluding remarks}
Finally, it should be noted that the rather simple argument used in this paper is in contrast with the much more sophisticated method used by B. Kats \cite{Ka} in solving Riemann problems on fractal curves (see also \cite{ABK,AB}). In such papers an appropriate domain integral replaced the role of the traditional Cauchy transform.  Although the Whitney extension is a crucial ingredient in both methods, the explicit use of the constructive algorithm of the snowflakes allows us to simplified the arguments as well as the calculations. The idea of the present paper is more in the direction of the results of Harrison and Norton  \cite{HN} where an integration theory of differential forms is developed by using  an approximating sequence of polyhedral $n$-chains.    

Of course we are aware that our achievement, in principle, is limited to Koch-type fractals. But at the same time we are of the opinion that the idea of the method can be used in the much more general case of the so-called deterministic fractals \cite{Br}.

\end{document}